\documentclass{aptpub}
\usepackage{amsmath,amssymb,epsfig,xspace}

\authornames{H.-O.\ GEORGII AND E.\ BAAKE}
\shorttitle{Multitype branching: the ancestral types}

\newcommand{\Definition}{\bigskip\noindent{\bfseries Definition: }}

\numberwithin{equation}{section}

\def\ba{\begin{array}}
\def\ea{\end{array}}
\def\be{\begin{equation} \label}
\def\ee{\end{equation}}
\def\bea{\begin{eqnarray*}}
\def\eea{\end{eqnarray*}}
\def\beal{\begin{eqnarray} \label}
\def\eeal{\end{eqnarray}}
\def\bit{\begin{itemize}}
\def\eit{\end{itemize}}
\def\ben{\begin{enumerate}}
\def\een{\end{enumerate}}

\def\inv{^{-1}}
\def\ti{\to\infty}
\def\Surv{\Omega_{\mathrm{surv}}}
\newcommand{\tper}[1][t]{#1,\mathrm{per}}
\newcommand{\toti}[1][n]{\;\mathrel{\mathop{\longrightarrow}\limits_{#1\ti}}\;}
\renewcommand{\sp}[1]{\langle #1\rangle}

\def\lo{\mathopen{]}\/}  %links offen
\def\ro{\/\mathclose{[}} %rechts offen

\def\RR{\mathbb{R}}
\def\NN{\mathbb{N}}

\def\EE{\mathbb{E}}
\def\FF{\mathbb{F}}

\def\ZZ{\mathbb{Z}}
\def\PP{\mathbb{P}}
\def\XX{\mathbb{X}}

\def\a{\alpha}
\def\d{\delta}
\def\eps{\varepsilon}
\def\th{\vartheta}
\def\k{\kappa}
\def\ph{\varphi}
\def\s{\sigma}
\def\t{\tau}
\def\Gam{\Gamma}
\def\S{\mathrm{\Sigma}}

\def\bmu{\boldsymbol{\mu}}
\def\bnu{\boldsymbol{\nu}}

\newcommand\bd[1]{{#1}}

\def\0{\mathbf{0}}
\def\L{\bd{L}}

\def\1{\mathsf{1}}
\def\p{\mathsf{p}}

\def\A{\mathsf{A}}
\def\M{\mathsf{M}}
\def\G{\mathsf{G}}

\def\hT{\widehat{T}}
\def\hW{\widehat{W}}
\def\hX{\widehat{X}}
\def\hN{\widehat{N}}
\def\hEE{\widehat{\EE}}
\def\hPP{\widehat{\PP}}

\def\hp{\widehat{\p}}
\def\htau{\widehat{\t}}

\def\tp{\widetilde{\p}}
\def\tx{\tilde{x}}

\def\l{\ell}
\def\cP{\mathcal{P}}
\def\cF{\mathcal{F}}
\def\cT{\mathcal{T}}
\def\fP{\mathfrak{P}}

\begin{document}

\bibliographystyle{apt.bst}

%\gii{Titel gek"urzt:  in continuous time}

\title{Supercritical multitype branching processes:\\
the ancestral types of typical individuals}

\vspace{8mm}

\authorone[Universit\"at M\"unchen]{Hans-Otto Georgii}
\addressone{Institut f\"ur Mathematik, Universit\"at M\"unchen,
            Theresienstr.\ 39,  D-80333 M\"unchen, Germany, 
            email: georgii@mathematik.uni-muenchen.de} 
\authortwo[Universit\"at Greifswald]{Ellen Baake}
\addresstwo{Institut f\"ur Mathematik und Informatik, Universit\"at
            Greifswald, Jahnstr.\ 15a, D-17487 Greifswald, Germany, 
            email: ellen.baake@uni-greifswald.de}

\begin{abstract}
For supercritical multitype Markov
branching processes in continuous time, we 
investigate the evolution of types along those lineages that survive up to 
some time $t$. We establish almost-sure convergence theorems for both time 
and population averages of ancestral types (conditioned on non-extinction),
and  identify the mutation process describing the type evolution along
typical lineages. An important tool is a representation of the family tree 
in terms of a suitable size-biased tree with trunk. As a by-product, this 
representation allows a `conceptual proof' (in the sense of \cite{KLPP97}) 
of the continuous-time version of the Kesten-Stigum theorem.
\end{abstract}



\section{Introduction}
\label{sec:intro}

Looking at the time evolution of a population one has two possible
perspectives:
either forward or backward in time. In the first case one observes the
characteristics of the population at a given time $t$ and asks for its
behaviour as $t$ increases to infinity. A classical model
that describes the unrestricted reproduction of independent individuals
is the (multitype) branching process, and a principal
result in the supercritical case is the Kesten-Stigum theorem
\cite{KeSt66}, which  describes the population size and relative frequencies
of types; see Theorem \ref{KS} for the precise statement.
A different situation arises if the population size is kept constant;
this leads to certain interacting particle systems,
like the Moran model
and its relatives (for review, see \cite{Durr02}).
By way of contrast, the backwards -- or retrospective -- aspect of the
population concerns the lineages extending back into the past
from the presently
living individuals  and asks for the characteristics of the
ancestors along such lineages. One famous example is
Kingman's coalescent (see \cite{King82a,King82b}, and \cite{Moeh00} for
a review),  the backward version of the Moran model.
As was observed e.g.\ by Jagers
\cite{Jage92} and Jagers and Nerman \cite{JaNe96}, it is also rewarding
to study the backward aspects of multitype branching processes;
this point of view has turned out as crucial in recent biological
applications \cite{HRWB02}.
It is the aim of this article to pursue this last line of research further.
We do so in continuous time because this gives us the
opportunity to  transfer  some powerful methods recently 
developed for discrete time.
We also concentrate on the supercritical case.

Specifically, we consider the individuals alive at some time
$t$ and investigate the types of their ancestors at an earlier time, $t-u$.
We will show the following.
\bit
\item When $t$ resp. $t$ and $u$ tend to infinity,
both time average and population average
of ancestral types converge to a particular distribution $\alpha$
almost surely on non-extinction (Theorems \ref{pop_average} and
\ref{time_average}).
\eit
This $\alpha$ will be called the  \emph{ancestral distribution of types}; its
components are $\alpha_i = \pi_i h_i$, where $\pi$ and
$h$ are the (properly normalized) left and right Perron-Frobenius
eigenvectors of the generator of the first-moment matrix.

More detailed information about the evolution of types
along ancestral lineages is obtained through what we would like to call the
\emph{retrospective mutation chain}, a particular continuous-time
Markov chain on the type space with $\alpha$ as its invariant distribution.
We will show:
\bit
\item For all individuals alive at time $t$ up to an asymptotically
negligible fraction, the time averaged empirical type evolution process
tends in distribution to the stationary retrospective mutation chain,
in the limit as $t\ti$, almost surely on non-extinction (Theorem
\ref{mutationhistory}).
\eit

One basic ingredient of our reasoning is a law of large numbers for
population averages; see Proposition \ref{dLLN}. A second crucial ingredient
is a representation of the family tree in terms of a size-biased tree with
trunk (with the retrospective
mutation chain running along the trunk); see Theorem
\ref{size_bias}. This representation
%, which allows to deduce properties of the whole population from those of
%the trunk alone,
is the continuous-time analogue of the size-biased tree representation
introduced by Lyons, Pemantle and Peres \cite{LPP95} and Kurtz, Lyons,
Pemantle and Peres \cite{KLPP97}. In passing, it allows us to extend their
conceptual proof of the Kesten-Stigum theorem to continuous time.
The third ingredient is the Donsker-Varadhan large deviation principle
for the retrospective mutation chain \cite{DoVa75a,DoVa83}.
This implies a large deviation
principle for the typical type evolution along the surviving lineages in the
tree  -- see Theorem \ref{LDP}.

This paper is organized as follows. In the next section we recall the
construction of the family tree for multitype branching processes in
continuous time. Section \ref{sec:results} contains the precise statement
of results. Section \ref{sec:bias_tree} is devoted to the size-biased tree
with trunk, and the proofs of the main results are collected in Section
\ref{sec:proofs}.

\section{The branching process and basic facts}
\label{sec:model}
\setcounter{equation}{0}

We consider a continuous-time multitype branching process as described
in Athreya and Ney \cite[Ch.~V.7]{AtNe72}. To fix the notation we recall
the basic setting here.

Let $S$ be a finite set of types. An individual of type $i \in S$
lives for an exponential time with parameter $a_i>0$,
and then splits into a  random offspring $N_i = (N_{ij})_{j \in S}$
with distribution $\p_i$ on $\ZZ^{S}_+$ and finite means
$m_{ij}:= \EE(N_{ij})$ for all $i,j \in S$; here, $N_{ij}$ is the
number of $j$-children, and $\ZZ_+=\{0,1,\ldots\}$.
We assume that the  mean offspring matrix
$\M=(m_{ij})_{i,j \in S}$ is irreducible.

According to Harris \cite[Ch.~VI]{Harris63}, the associated random
family tree can be constructed as follows. Let
$\XX = \bigcup_{n \geq 0} \XX_n$, where $\XX_n$ describes
the virtual $n$'th generation.  That is, $\XX_0=S$, and $i_0\in \XX_0$
specifies the type of the root, i.e., the founding ancestor. Next,
$\XX_1 =S \times \NN$, and the element $x=(i_1,\l_1) \in \XX_1$ is the
$\l_1$'th $i_1$-child of the root. Finally, for $n>1$,
$\XX_n =S^n \times \NN^n$, and $x = (i_1,\ldots,i_n;\l_1,\ldots,\l_n)
\in \XX_n$ is the $\l_n$-th $i_n$-child of its parent
$\tx= (i_1, \ldots, i_{n-1};\l_1, \ldots, \l_{n-1})$; see Fig.~1.
We write $\s(x)=i_n$ for the type of $x \in\XX_n$.
With each $x \in \XX$ we associate
\bit
\item its random life time $\t_{x}$, distributed exponentially
with parameter $a_{\s(x)}$, and

\item its random offspring $N_{x} = (N_{x,j})_{j \in S} \in \ZZ_+^{S}$
with distribution $\p_{\s(x)}$
\eit
such that the family $\{ \t_{x}, N_{x} : \, x \in \XX \}$
is independent.

The random variables $N_{x}$ indicate which of the virtual
individuals $x \in \XX$ are actually realized, namely those in the random set
$X = \bigcup_{n \geq 0} X_n$ defined recursively by
\[
X_0 = \{i_0\}, \quad X_n = \{x = (\tx; i_n, \l_n) \in
\XX_n: \; \tx \in X_{n-1}, \,\l_n \leq N_{\tx,i_n} \},
\]
where $i_0$ is the prescribed type of the root.
The random variables $\t_{x}$ provide the proper time scale.  Namely,
for $x \in X$, let the splitting times $T_{x}$ be defined recursively
by $T_{x} = T_{\tx} + \t_x$ with $T_{\tilde{i_0}} := 0$.  The
lifetime interval of $x \in X$ is then
$[T_{\tx},T_{x}\ro$.  Hence $X(t) = \{x \in X: \,T_{\tx}\le t < T_{x}\}$
is the population at time $t$.  One may visualize the
resulting tree by identifying each $x\in X$ with an edge from
$\tx$ to $x$ with length $\t_{x}$ in the direction of time.

\begin{figure}[htb]
\centerline{\epsfig{file=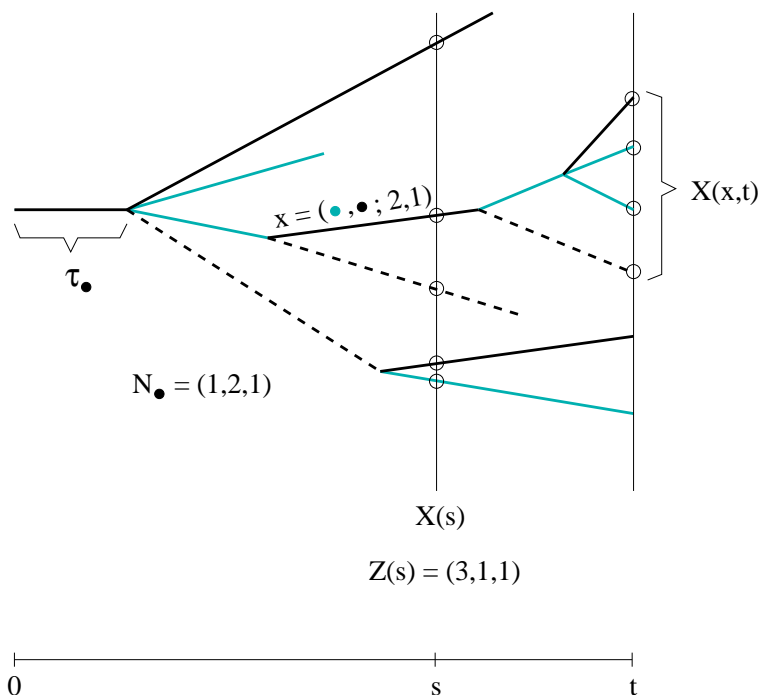,width=100mm}}
\caption{
A realization of the branching process. Types are indicated
by different line types, 
indexed in the order (black, grey, dashed), counted from top to bottom,
and symbolized by filled circles.
The set $X(s)$ consists of all edges that intersect the vertical line at $s$;
the set $X(x,t)$ consists of all edges that emanate from $x$ and
hit the vertical line at $t$. $Z(s)$ counts the type frequencies
in the population $X(s)$.}
\end{figure}

The family tree is completely determined by the process
$X[0,\infty\ro:=\big ( X(t) \big )_{t \geq 0}$ which is a random element of
$\Omega:= D\big ( [0,\infty \ro,\fP_f(\XX)\big )$, the Skorohod space of all
c\`{a}dl\`{a}g functions on $[0,\infty\ro$ taking values in the
(countable) set $\fP_f(\XX)$ of all
finite subsets of $\XX$.  We write $\PP^{i}$ for the distribution of
$X[0,\infty\ro$ on $\Omega$ when the type of the root is $i_0=i$, and
$\EE^i$ for the associated expectation. If $i_0$
is chosen randomly with distribution $\nu$, we write $\PP^\nu$ and
$\EE^\nu$. We will often
identify $X[0,\infty\ro$ with the canonical process on
$\Omega$.

For $0<s<t$ and $y\in X(t)$ we write $y(s)$ for its unique ancestor living at
time $s$. On the other hand, for $x\in X(s)$ we let
\be{Xxt}
X(x,t)=\big\{ y\in\XX: xy\in X(t)\big\}\
\ee
denote the set of descendants of $x$ living at time $t$;
cf.\ Fig.~1.
In the above,
the concatenation $xy$ of two strings $x,y\in\XX$ is defined in
the obvious way, and the empty string is considered as an ancestor
of type $\s(x)$; i.e., $X(x,t)=\{\s(x)\}$ as long as $x\in X(t)$.
By the loss-of-memory property of the exponential distributions,
the descendant
trees $X(x,[s,\infty\ro)=(X(x,t))_{t\ge s}$ with $x\in X(s)$ are conditionally
independent given $X[0,s]$,  with distribution $\PP^{\s(x)}$.
We will also consider the counting measures
\be{Zxt}
Z(t)= \sum_{x\in X(t)} \delta_{\s(x)}\,,\quad
Z(x,t)= \sum_{y\in X(x,t)} \delta_{\s(y)}
\ee
on $S$, where $\delta_i$ is the Dirac measure at $i$.
$Z(t)$ and $Z(x,t)$ count the type frequencies in the
population $X(t)$ resp. the subpopulation $X(x,t)$ of $x$-descendants.
In particular, $Z_j(t)$ is the cardinality of
$X_j(t) = \{ x \in X(t): \s(x) = j \}$, the subpopulation of
type $j\in S$, and  $\|Z(t)\|:=\sum_{j\in S}Z_j(t) = |X(t)|$ is the total
size of the population.

It is well-known (cf.\ \cite{AtNe72}, p.~202, Eq.~9) that
$\EE^i( Z_j(t) ) = (e^{t \A})_{ij}$ for all
$i,j \in S$,  where the generator matrix $\A=(a_{ij})_{i,j\in S}$ is
given by
\be{a}
  a_{ij} = a_{i} ( m_{ij} - \delta_{ij} )\,.
\ee
By the irreducibility of $\M$, $\A$ is also irreducible, so that the first
moment matrix $(\EE^i( Z_j(t) ) )_{i,j\in S}$ has positive
entries for any $t>0$. (This property is often
called `positive regularity', see \cite[p.~202]{AtNe72}.)
Perron-Frobenius theory then tells us that the matrix
$\A$ has a principal eigenvalue $\lambda$
(i.e., a real eigenvalue exceeding the real parts of all other eigenvalues),
and associated positive left and right
eigenvectors $\pi$ and $h$ which will be normalized s.t.\
$\sp{ \pi, \1 } = 1 = \sp{ \pi ,h }$.
Here we think of the row vector $\pi$
as a probability measure, of the column vectors $h$ and
$\1=(1,\ldots,1)^T$ as
functions on $S$, and of the scalar product $\sp{ \pi,h } =
\sum_i \pi_i h_i$ as the associated
expectation.
We are mainly interested in the supercritical case
$\lambda>0$. In this case we write
\[
\Surv :=\{ X(t)\ne\varnothing \mbox{ for all }t>0\}
\]
for the event that the population survives for all times.

It is a remarkable fact that the almost-sure behaviour of the family
tree is, to a large extent, already determined by the the global
quantities $\lambda, \pi, h$. One prominent example is the
following continuous-time version of the Kesten-Stigum theorem
(see \cite{KeSt66} for the discrete-time original,
\cite{Athr68} for the continuous-time
version, and \cite{KLPP97} for the recent discrete-time conceptual proof).

\begin{thm}[Kesten-Stigum]\label{KS} Consider the supercritical case
$\lambda>0$.

{\rm(a)} For all $i\in S$ we have
\[
\frac{1}{ |X(t)|}\sum_{x\in X(t)} \delta_{\s(x)}
= \frac{Z(t)}{\|Z(t)\|} \toti[t] \pi
\quad\mbox{$\PP^i$-almost surely on $\Surv$.}
\]

{\rm(b)} There is a nonnegative random variable $W$ such that
\[
\lim_{t\ti} Z(t)\, e^{-\lambda t} = W\,\pi
\quad\mbox{$\PP^i$-almost surely for any $i\in S$}\,,
\]
and $\PP^i(W>0)>0$ for all $i$ if and only if
\be{ZlogZ}
\EE(N_{ij}\log N_{ij})<\infty\quad\mbox{for all $i,j\in S$.}
\ee
In this case, $\{W>0\} = \Surv$ $\PP^i$-almost surely, and
$h_i=\EE^i(W)$.
\end{thm}

For the sake of reference we provide here a full proof extending
the conceptual discrete-time proof of \cite{KLPP97} to our
continuous-time setting.
Assertion (a) %(to be proved in Section \ref{subsec:LLN})
reveals that the left eigenvector $\pi$
holds the asymptotic proportions of the types in the population,
and statement (b) %(which will be proved in Section \ref{subsec:treetrans})
implies that
\[
\lambda = \lim_{t\ti} \frac1t \log |X(t)|
\]
is the almost sure exponential growth rate of the population
in the case of survival. In fact, this statement does not require condition
(\ref{ZlogZ}); see the proof of Theorem \ref{mutationhistory} in Section
\ref{subsec:LDP}. The $i$-th
coordinate $h_i$ of the right eigenvector $h$ measures the
long-term fertility of an $i$-individual.
In fact, $h_i$ is also characterized
by the limiting relation
\be{popsize}
\EE^{i} \big(|X(t)|\big)\, e^{-\lambda t}
\to h_{i}\quad\mbox{as $t \ti$;}
\ee
cf. Remark \ref{conseq}(a) below.

\section{Results}
\label{sec:results}
\setcounter{equation}{0}

We still consider the supercritical case $\lambda>0$.  We are
interested in the mutation behaviour of the population tree.  More
specifically, we ask for the behaviour of the sequence of types
along a typical branch of this tree.  It turns out that this behaviour
is again completely determined by the global quantities $\lambda,
\pi, h$.  A key role is played by the probability vector
$\a=(\a_i)_{i\in S}$ with components $\a_i = \pi_i h_i$.
As observed by Jagers \cite[Corollary 1]{Jage92}, Jagers and
Nerman \cite[Prop.~1]{JaNe96}, and Hermisson et al.\ \cite{HRWB02},
this probability vector describes the distribution of
ancestral types of an equilibrium population with type frequencies
given by $\pi$.  The vector $\a$ will therefore be called the
\emph{ancestral distribution}.  Our results below shed some additional
light on the significance of $\a$.

To begin, we consider a typical individual $x\in X(t)$ alive at some large time
$t$ and ask for the type $\s(x(t{-}u))$ of its ancestor $x(t{-}u)$ living at
some earlier time $t{-}u$. We find that $\s(x(t{-}u))$ is asymptotically
distributed according to $\a$.
Specifically, let $0<u<t$ and
\be{eq:ancestral-average}
A^u(t) =\frac1{|X(t)|}\sum_{x\in X(t)}\delta_{\s(x(t{-}u))}
\ee
be the empirical ancestral type distribution at time $t{-}u$ taken over
the population $X(t)$. (Of course, this definition requires that
$X(t)\ne\varnothing$.)

\begin{thm}[Population average of ancestral types]\label{pop_average}
Let $\lambda>0$ and $i\in S$.  Then
\be{eq:pop_average}
\lim_{u\ti}\lim_{t\ti} A^u(t)
%\frac1{|X(t)|}\sum_{x\in X(t)}\delta_{\s(x(t-u))}
=\a
\quad\mbox{$\PP^i$-almost surely on $\Surv$.}
\ee
\end{thm}

The proof will be given in Section \ref{subsec:LLN}. We would like
to remark that a slightly weaker result under slightly stronger
conditions (convergence in probability under assumption (\ref{ZlogZ}))
follows immediately from Corollary~4 of Jagers and Nerman \cite{JaNe96},
where  very general population averages are considered.

\begin{rem}\label{rem:pop_average}\rm
Assertion (\ref{eq:pop_average}) means that, for each $j\in S$,
the average
\[
A^u_j(t) = \frac1{|X(t)|}\sum_{x\in X(t)} I\{\s(x(t{-}u))=j\}
\]
(with $I\{.\}$ denoting the indicator function)
converges to $\a_j$ $\PP^i$-almost surely on $\Surv$ as
$t\ti$ and $u\ti$ in this order. Letting $s=t{-}u$, we can rewrite
this average in the form
\[
\sum_{x\in X_j(s)} |X(x,t)| \Big/
\sum_{x\in X(s)} |X(x,t)|\,,
\]
where $X(x,t)$ is given by (\ref{Xxt}).
The numbers $|X(x,t)|$ with $x\in X_j(s)$ are
i.i.d. with mean $\EE^j(|X(u)|)$. Assuming the validity of a
law of large numbers and using
Theorem \ref{KS}(a) and Eq.~(\ref{popsize}), we can conclude that
the average above
converges to $\pi_j h_j/\sp{ \pi,h}=\a_j$ as
$s, u\ti$. This explains the particular structure of the ancestral
distribution $\a$.
\end{rem}

In our next theorem we ask for the time average of types along the line
of descent leading to a typical $x\in X(t)$. This time average is given
by the empirical distribution
\[
\L^x(t)=\frac{1}{t} \int_0^t \delta_{\s(x(s))}\, ds
\]
of the process $\s(x[0,t])=\big(\s(x(s))\big)_{0\le s\le t}$.
Note that $\L^x(t)$ belongs to the simplex $\cP(S)$ of all probability
vectors on $S$; $\cP(S)$ will be equipped with the usual total
variation distance $\|\cdot\|$.  To describe the behaviour of
$\L^x(t)$ for a typical $x\in X(t)$ we have to step one level
higher and to consider the empirical distribution of $\L^x(t)$
taken over the population $x\in X(t)$.  This empirical distribution
belongs to $\cP(\cP(S))$, the set of probability measures on $\cP(S)$,
which will be equipped with the weak topology.

\begin{thm}[Time average of ancestral types]\label{time_average}
Let $\lambda>0$ and $i\in S$. Then
\be{eq:time_average}
\lim_{t\ti} \frac1{|X(t)|}\sum_{x\in X(t)}\delta_{\L^x(t)}
=\delta_{\a} \quad\mbox{$\PP^i$-almost surely on $\Surv$.}
\ee
\end{thm}

\begin{rem}\label{rem:time_average}\rm
(a) According to the portmanteau theorem
\cite[p.~108, Th.~3.1]{EtKu86},
statement (\ref{eq:time_average}) is equivalent to the assertion that
\[
\lim_{t\ti} \frac1{|X(t)|}\sum_{x\in X(t)}I\{\,\L^x(t)\in F\}=0
\quad\mbox{for each closed $F\subset\cP(S)$ with $\a\notin F$}
\]
$\PP^i$-almost surely on $\Surv$, and it is sufficient to check this
in the case when $F=\{\nu\in\cP(S):\|\nu-\a\|\ge \eps\}$
with arbitrary $\eps>0$. The theorem therefore asserts that,
for all individuals $x\in X(t)$ up to an asymptotically negligible
fraction, the ancestral type average $\L^x(t)$ is close to $\a$.

(b) Theorem \ref{time_average} involves a population average of time averages.
So one may ask whether the averaging of population and time can be interchanged.
It follows from Theorem \ref{pop_average} that this is indeed the case:
\[
\lim_{t\ti} \frac1t \int_0^t  \delta_{A^u(t)}\; du
=\delta_{\a}
\]
almost surely on $\Surv$.
\end{rem}

Theorem \ref{time_average} is in fact a corollary of our next theorem
which considers the complete mutation history along a typical line
of descent. To state this result we need some preparations.
We introduce first the mutation process on $S$ which will turn out to
describe the
time-averaged mutation behaviour along an ancestral line.

\Definition
The \emph{retrospective mutation chain} is the
Markov chain $ (\s(t))_{t\geq0}$ on $S$ which stays in a state
$i\in S$ for an exponential holding time with parameter $a_i{+}\lambda$ and
then jumps to $j\in S$ with probability
\[
    p_{ij}  =  \frac{m_{ij}\, h_j}{(1{+}\lambda/a_i) h_i} \,.
\]
That is, the generator $\G = (g_{ij})_{i,j \in S}$ of $ (\s(t))_{t\geq0}$
is given by
\[
   g_{ij}  = (a_i{+}\lambda) (p_{ij} - \delta_{ij}) =
   h_i\inv (a_{ij} - \lambda \d_{ij}) h_j.
\]

\bigskip
We note that $\G$ is indeed a generator because
$a_i\sum_{j \in S} m_{ij} h_j =  \sum_{j \in S}
   (a_i \delta_{ij}+a_{ij}) h_j  = (a_i{+}\lambda) h_i $
by (\ref{a}).
Since $\M$ is irreducible by assumption, $\G$ is irreducible as well.
It is also immediate that the ancestral distribution $\a$ is
the (unique) stationary distribution of $\G$.
The retrospective mutation chain was  identified by Jagers
\cite[p.~195]{Jage89} and may be interpreted as the
forward version of the backward Markov chain 
\cite[Proposition 1]{JaNe96} that results from
picking individuals randomly from the stationary
type distribution $\pi$ and following their lines of descent
backward in time. This gives the transition rates
\be{tG}
  \bar g_{ij} = \pi_j (a_{ji} - \lambda \d_{ij}) \pi_i^{-1}
  = \a_j g_{ji} \a_i^{-1}\,,
\ee
which corresponds to
the time reversal of the retrospective mutation chain.

To set up the stage for Theorem \ref{mutationhistory} we let
$\S=D(\RR,S)$ denote
the space of all doubly infinite c\`adl\`ag paths in $S$.
$\S$ will be equipped with the usual Skorohod topology which turns
$\S$ into a Polish space; see e.g.\ \cite{EtKu86},
Section 3.5 and in particular Th.~5.6, for the case of the time interval
$[0,\infty\ro$. The associated Borel $\s$-algebra coincides with the
$\s$-algebra generated by the evaluation maps $\S\ni\s\to \s(t)$,
$t\in\RR$ \cite[p.~127, Prop.~7.1]{EtKu86}. The time shift
$\th_s$ on $\S$ is defined by
\[
\th_s \s(t) = \s(t+s)\,,\qquad s,t\in\RR,\,\s\in\S.
\]
We write $\cP_\Theta(\S)$ for the set of all probability measures
on $\S$ which are invariant under the shift group $\Theta=(\th_s)_{s\in\RR}$.
Endowed with the weak topology, $\cP_\Theta(\S)$ is a Polish space
\cite[p.~101, Th.~1.7]{EtKu86}.

Next we introduce the time-averaged type evolution process of an individual
in the population tree. For $t>0$ and $x\in X(t)$ we let $\s(x)_{\tper}\in\S$
be defined by
\be{sigmaper}
\s(x)_{\tper}(s) = \s(x(s_t))\,,\quad s\in\RR\,,
\ee
where $s_t$ is the unique number in $[0,t\ro$ with $s\equiv s_t \mod t$.
That is, $\s(x)_{\tper}\in\S$ is the periodically continued type history
of $x$ up to time $t$. The time-averaged type evolution of $x$ is then
described by the \emph{empirical type evolution process}
\be{Rt}
R^x(t)= \frac 1t \int_0^t \delta_{\th_s\s(x)_{\tper}}\,ds
\ \in \cP_\Theta(\S).
\ee

We are interested in the typical behaviour of $R^x(t)$ when
$x$ is picked at random from $X(t)$, the population at time $t$.
This is captured in their empirical distribution, i.e., the
population average
\be{Ht}
\Gam(t) :=\frac1{|X(t)|}\sum_{x\in X(t)}\delta_{R^x(t)}\,.
\ee
(As before, this definition requires that $X(t)\ne\varnothing$.)
$\Gam(t)$ is a random element of $\cP(\cP_\Theta(\S))$, the set of
all probability measures on the Polish space $\cP_\Theta(\S)$,
which is again equipped with the weak topology. In
Section \ref{subsec:LDP} we will prove:

%Typical ancestral mutation histories
\begin{thm}[Typical ancestral type evolution]\label{mutationhistory}
Let $\lambda>0$ and $i\in S$. Then
\be{eq:mutationhistory}
\lim_{t\ti} \Gam(t)
=\delta_{\bmu} \quad\mbox{$\PP^i$-almost surely on $\Surv$,}
\ee
where $\bmu\in\cP_\Theta(\S)$ is the distribution of the stationary
(doubly infinite) retrospective mutation chain $(\s(t))_{t\in\RR}$ with generator $\G$ and invariant distribution $\a$.
\end{thm}

\begin{rem}\label{rem:mutationhistory}\rm
As in Remark \ref{rem:time_average}(a), the portmanteau theorem implies
that (\ref{eq:mutationhistory})
is equivalent to the assertion that, $\PP^i$-almost surely on $\Surv$,
$\Gam(t)(F)\to0$ for every closed $F\subset\cP_\Theta(\S)$ such that
$\bmu\notin F$. Writing $d(\cdot,\cdot)$ for any metric metrizing
the weak topology on $\cP_\Theta(\S)$ this in turn means that, for each $\eps>0$,
\[
\lim_{t\ti} \frac1{|X(t)|}\sum_{x\in X(t)}I\big\{d(R^x(t),\bmu)\ge\eps\big\}=0
\]
$\PP^i$-almost surely on $\Surv$. The theorem therefore states that,
for all individuals $x\in X(t)$ up to an asymptotically negligible fraction,
the time-averaged ancestral type evolution process $R^x(t)$ is close
to $\bmu$ in the weak topology.
Theorem \ref{mutationhistory} also highlights the restrospective
nature of our mutation chain:
it describes the evolution of types along those lines of descent which
survive until time $t$ (and thus can be seen when a time-$t$ individual
looks back into the  past).
\end{rem}

\section{Size-biasing of the family tree}
\label{sec:bias_tree}
\setcounter{equation}{0}

In this section we construct a continuous-time version of the size-biased
multitype Galton-Watson tree as introduced by Lyons, Pemantle, Peres, 
and Kurtz \cite{LPP95,KLPP97}. Informally, this is a  tree
with a randomly selected trunk (or spine)
along which time runs at a different rate and offspring is weighted
according to its size; in particular,
there is always at least one offspring along the trunk so that the trunk survives forever.
The children off the trunk get ordinary (unbiased) descendant trees (the  bushes).
It will turn out that the trunk of the size-biased tree describes the evolution along
a typical ancestral line that survives up to some fixed time.
The construction is not confined to the supercritical case; that is,
in this section $\lambda$ can have arbitrary sign.

First of all, for each type $i\in S$ we introduce the size-biased
offspring distribution
\be{sbd}
  \hp_i(\k) = \frac{\sp{\k, h}\, \p_i(\k)}
                  {c_i \,h_i }\,,\qquad \k\in\ZZ_+^S,
\ee
where $\sp{ \k, h } = \sum_j \k_j h_j$ and $c_i=1+ \lambda/a_i$ is a
normalizing constant. $\hp_i$ will serve as the offspring distribution
of an $i$-individual on the trunk; it
is indeed a probability distribution since
\[
\sum_{\k\in\ZZ_+^S} \sp{\k, h}\, \p_i(\k) = \sum_{j\in S} m_{ij} h_j
= \sum_{j \in S}
   (\delta_{ij}+a_{ij}/a_i) h_j  = c_i \,h_i
\]
by (\ref{a}); note that $c_i$ is automatically positive.
Next, when an $i$-individual on the trunk has offspring
$\hN_i=(\hN_{ij})_{j\in S}$ with distribution $\hp_i$, 
one of these offspring is chosen as the successor on the trunk,
where children are picked with probability proportional to $h_j$
when their type is $j$. That is, the successor is of type $j$
with probability $\hN_{ij}\; h_j / \sp{\hN_i,h}$ for a given
offspring, and with probability
\[
p_{ij}= \EE\Big(\frac{\hN_{ij}\; h_j}{\sp{ \hN_i,h}}\Big)=
\frac{m_{ij}\,h_j}{ c_i\,h_i}
\]
on average.
\emph{These are precisely the jump probabilities of the retrospective mutation
chain.} Finally, the lifetime of an $i$-individual on the trunk will be exponential
with parameter $a_i{+}\lambda$, which coincides again with the holding time
parameter of the retrospective mutation chain. A corresponding 
embedded chain combined with size-biased waiting times also occurs when
more general non-Markovian populations (i.e., with waiting times
deviating from the exponential distribution) are traced backwards,
see \cite[Proposition 1]{JaNe96}.

We now construct the size-biased tree in detail.
Let $\{\t_{x}, N_{x} : \, x \in \XX \}$ be as in Section~\ref{sec:model}
and, independently
of this, a sequence $\{ \htau_n, \hN_n, \xi_n: n \geq 0 \}$ of
random variables with values in $\lo 0, \infty \ro \,, \ZZ_+^{S}, \XX$
respectively such that, for a given type $i_0=i$ of the
root,   $\xi_0=i$ and
\begin{itemize}
\item $\htau_0, \hN_0$ are independent, $\htau_0$ has exponential
distribution with parameter $a_{i} {+}\lambda$,  $\hN_0$ has distribution
$\hp_{i}$, and $\xi_1$ has conditional distribution
\[
   P\big(\xi_1 = (i_1,\l_1) | \hN_0,\htau_0\big)
   = \frac{h_{i_1}}{\sp{ \hN_0 , h}}
    \, I{\{\l_1 \leq \hN_{0,i_1}\}}
\]
for all $(i_1,\l_1)\in\XX_1$.

\item For any $n \geq 1$, conditionally on $\FF_{n-1} =
\s\{\htau_k,\,\hN_k,\,\xi_{k+1}: \, k < n \}$, $\htau_n,\hN_n$ are independent
and follow an exponential law with parameter $a_{\s(\xi_n)} {+}\lambda$
resp.\ the law $\hp_{\s(\xi_n)}$,
and
\[
  P \big ( \xi_{n+1} = (\xi_n;i_{n+1},\l_{n+1}) | \FF_{n-1},\htau_n,\hN_n \big )
  = \frac{h_{i_{n+1}}}{\sp{ \hN_n , h }}  I{\{\l_{n+1}
  \leq \hN_{n,i_{n+1}}\}}
\]
for all $(i_{n+1},\l_{n+1})\in S\times\NN$, i.e., $\xi_{n+1}$ is a
child of $\xi_n$ selected randomly with weight
proportional to $h_{\s(\xi_{n+1})}$.
\end{itemize}
Define $\hX = \bigcup_{n \geq 0} \hX_n \subset \XX$ recursively by
$\hX_0 = \{i\}$ and $\hX_n = \hX_n^{\sharp} \cup \hX_n^{\flat}$ with
\[
 \hX_n^{\sharp} =
    \{ (\xi_{n-1}; i_n,\l_n) \in \XX_n: \l_n \leq \hN_{n-1,i_n} \},
\]
the offspring of $\xi_{n-1}$, and
\[
 \hX_n^{\flat} = \{ (\tx; i_n,\l_n) \in \XX_n:  \tx \in \hX_{n-1}
 \setminus \{\xi_{n-1}\}, \l_n \leq N_{\tx,i_n} \}
\]
the offspring of all other individuals in $\hX_{n-1}$. (Note that in the last
display there is no hat on $N$; that is, the bushes have unbiased offspring.)
The split times $\hT_{x}$ are given by
$\hT_{\xi_0} = \htau_0$, $\hT_{\xi_n} = \hT_{\xi_{n-1}}
+ \htau_n$ for $n \geq 1$, and $\hT_{x} = \hT_{\tx} + \t_{x}$ if
$x \in \hX \setminus \{\xi_n: \, n \geq 0 \}$. (Again, in the latter case there is no hat on $\t$, meaning that the individuals off the trunk have unbiased life times.) The total population at time $t$ is then given by
\[
   \hX(t) = \{x \in \hX: \, \hT_{\tx} \leq t < \hT_{x} \}\,.
\]
The selected trunk individual at time $t$ is
$\xi(t) = \xi_n$ if $\hT_{\xi_{n-1}} \leq t < \hT_{\xi_n}$, and the process
$\big(\hX(t),\xi(t) \big )_{t \geq 0}$ in $\Omega_*:=
D\big([0,\infty\ro,\fP_f(\XX)\times \XX \big )
=\Omega\times D\big([0,\infty\ro,\XX \big )$ describes the size-biased
tree with trunk $\big ( \xi(t) \big )_{t \geq 0}$.  As we have emphasized
above, the type
process along the trunk, $\s(t) := \s \big ( \xi(t) \big )$, is a copy
of the retrospective mutation chain as defined in Section \ref{sec:results}.
In contrast, the individuals off the trunk may be understood as a
branching process with immigration.
%where immigration occurs with
%rate $\lambda$, the excess branching rate of the trunk (which does not
%depend on the type).
%more precisely, during the time intervals for which the trunk has type $i$,
%with rate $a_i+\lambda$ one has immigration of independent clans with
%size distribution
%\[
%\hp_i'(\k) = \sum_{j\in S}\frac{(\k_j{+}1)h_j}{\sp{\k{+}\d_j,h}}
%\hp_i(\k{+}\d_j)\,,
%\quad \k\in\ZZ_+^S.
%\]

We write $\hPP_*^{i}$ for the distribution of
$\big ( \hX(t),\xi(t) \big )_{t \geq 0}$ on $\Omega_*$,
and $\hPP^{i}$ for its marginal, the distribution of
$\big ( \hX(t) \big )_{t \geq 0}$ on $\Omega$.
The representation theorem below establishes the relationship between 
$\PP^{i}$,
$\hPP_*^{i}$ and the retrospective mutation chain. We use the
shorthand $y[0,t]$ for a path $\big ( y(s) \big )_{0 \leq s \leq t}$.
 
\begin{thm}\label{size_bias}
Let $t>0$, $i \in S$, and 
$F:D\big ([0,t],\fP_f(\XX) \times \XX\big) \to [0,\infty\ro$ 
be any measurable function. Then one has
 \be{eq:size_bias}
  h_{i}\inv\, \EE^{i}
  \bigg (e^{-\lambda t}  \sum_{x \in X(t)} F\big(X[0,t],x[0,t]\big)\, h_{\s(x)} \bigg )
   = \hEE_*^{i} \Big ( F\big(\hX[0,t],\xi[0,t]\big) \Big ) \,.
 \ee
\end{thm}

\bigskip
Recall that this theorem is valid for arbitrary sign of $\lambda$.
The proof is postponed until Section \ref{subsec:treetrans}.
Here we discuss some immediate consequences and possible extensions.

%: rem:conseq
\begin{rem}\label{conseq}\rm
(a) Setting $F(X[0,t],x[0,t]) = I\big\{ \s  ( x(t)  ) = j \big\}\, h_j^{-1}$ in
(\ref{eq:size_bias}) and using the ergodic theorem for the retrospective
mutation chain $\s\big(\xi(t)\big)$ we obtain the Perron-Frobenius result
\be{PF}
  \EE^{i} \big ( Z_j(t) \big )\, e^{-\lambda t} =
  h_{i}\, \hPP_*^{i} \big ( \s(\xi(t)) = j \big )\, h_j^{-1}
  \toti[t] h_i\, \a_j\,h_j\inv
  = h_i\, \pi_j\,.
\ee
In particular, Eq. (\ref{popsize}) follows by summing over $j$.

(b) Taking any $F$ of the form $F(X[0,t],x[0,t]) = g(X[0,t])$
we conclude that
\[
h_{i}\inv\,  \EE^{i} \Big (W(t)\, g(X[0,t])  \Big )  =  \hEE^{i}\Big ( g(\hX[0,t]) \Big )
\]
with
\[
W(t) := \sp{ Z(t),h } \, e^{- \lambda t}\,.
\]
In particular, $h_{i} = \EE^{i} \big ( W (t) \big )$. Thus,
on the $\s$-algebra $\cF_t$ generated by $X[0,t]$,
$\hPP^{i}$ is absolutely continuous with respect to $\PP^{i}$
with density $W(t)/h_{i}$,
and $(W(t))_{t\ge 0}$ is a martingale with respect to $\PP^{i}$.
The latter statement is one of the standard facts of branching process theory;
see e.g. \cite{AtNe72}, p.~209, Theorem 1.

(c) Theorem (\ref{size_bias}) has the appearance of the Campbell theorem
of point process theory; see, e.g., \cite{MKM78}, pp.~14 \& 228.
To clarify the relation let $t>0$ be fixed and
\[
\Phi(t)= \big\{x[0,t]: x\in X(t) \big\}
\]
the finite random subset of $D\big([0,t],\XX\big)$ which describes the lineages that
survive until time $t$. Also, let $C^i_t$ be the measure on
$\fP_f\big(D\big([0,t],\XX\big)\big)\times D\big([0,t],\XX\big)$ with Radon-Nikodym
density $e^{\lambda t} \,h_i \,h_{\s(\xi(t))}\inv$ relative
to the joint distribution of
$\widehat{\Phi}(t)=\big\{x[0,t]: x\in \hX(t) \big\}$
and $\xi[0,t]$ under $\hPP^i_*$. Theorem (\ref{size_bias}) then implies that
\[
 \EE^{i}
 \Big ( \sum_{\psi \in \Phi(t)} F\big(\Phi(t),\psi\big) \Big )
= \int F(\Psi,\psi)  \, C^i_t(d\Psi,d\psi)
\]
for any measurable $F\ge 0$, i.e., $C^i_t$ is the Campbell measure of $\Phi(t)$
under $\PP^i$.
This assertion, however, is slightly weaker than Theorem (\ref{size_bias})
because $X[0,t]$ also includes the lineages that die out before time $t$.
\end{rem}

\begin{rem}\rm
In the above, the size-biased tree was constructed using the
right eigenvector $h$ as a weight on the types. As a matter of fact, the same
construction can be carried out when $h$ is replaced by an arbitrary weight
vector $\gamma\in\lo 0,\infty\ro^S$, and a representation theorem analogous to
Theorem (\ref{size_bias}) can be obtained. We discuss here only the special
case $\gamma\equiv1$ which is of particular interest,
and already appears in
\cite[Theorem 2]{GRW92} in the context of critical multitype
branching. The size-biased offspring
distribution associated with this case is
\[
  \tp_i(\k) = {\|\k\|\, \p_i(\k)}/{m_i }\,,\qquad \k\in\ZZ_+^S,
\]
where $\|\k\|=\sum_j\k_j$ is the total offspring and $m_i =\sum_j m_{ij}$
its expectation under $\p_i$. The lifetime of an
$i$-individual on the trunk is exponential with parameter $a_i m_i$, and
the successor on the trunk  is chosen among the children with
\emph{equal} probability.
Writing a tilde (instead of a hat) to characterize all quantities of the
associated size-biased tree, one arrives at the following counterpart of
(\ref{eq:size_bias}):
\be{eq:equal_size_bias}
\EE^{i}
 \bigg ( \sum_{x \in X(t)}  e^{- t\sp{ L^x(t),r}} \,
F\big(X[0,t],x[0,t]\big) \bigg ) =
\widetilde{\EE}_*^{i} \Big ( F\big(\widetilde{X}[0,t],\tilde{\xi}[0,t]\big) \Big ) \;.
 \ee
In the above, $r$ is the vector with $i$-coordinate $r_i=a_i(m_i-1)=\sum_j a_{ij}$,
the mean reproduction rate of type $i$. Accordingly, the expectation
$\sp{ L^x(t),r}$ is the \emph{mean reproduction rate along the lineage leading to $x$
at time $t$}. The type process along the trunk, $\tilde\s(t) := \s \big ( \tilde\xi(t) \big )$,
is the Markov chain with transition rates
$\tilde g_{ij} = a_i\, m_{ij} - m_i\,\d_{ij}$.
In view of the decomposition $ a_{ij}  = \tilde g_{ij} + r_i \delta_{ij}$,
this Markov chain
describes the pure mutation part of the type evolution.

On the left-hand side of (\ref{eq:equal_size_bias}), each individual
is weighted according to the mean fertility of its lineage. Indeed, suppose we are given a lineage
up to time $t$ of which we know only the intervals of time
spent in each state $i\in S$,
%but not the actual split times and offspring sizes.
%To determine the expected product of offspring sizes along the lineage,
and imagine that random split events
and independent random offspring sizes are distributed over $[0,t]$ with the
appropriate rates and distributions.
The number $\zeta_i$ of split events during the sojourn in state
$i$ is then Poisson
with parameter $a_i t \nu_i$, where $\nu_i$ is the fraction of time
spent in state $i$; and the expected total offspring at each of
these events is $m_i$. Since offspring sizes are independent,
the expected product of offspring sizes along the lineage then amounts
to $\prod_{i\in S}\EE \big (m_i^{\zeta_i}\big) =  e^{t \sp{\nu,r}}$.
A result similar to (\ref{eq:equal_size_bias}),  with an analogous
interpretation of the exponential factor, already appears in
\cite[p.~127]{CRW91} in the context of Palm trees for spatially
inhomogeneous branching. 

Here are some
consequences of (\ref{eq:equal_size_bias}):

(a) For $F(X[0,t],x[0,t]) = \exp\big[\,t\, \sp{\L^x(t),r}\big]\,I{\{\s(x(t)) = j\}} $,
Eq.~(\ref{eq:equal_size_bias}) becomes
\be{FC}
\EE^{i} \Big( Z_j(t) \Big)
%= \widetilde{\EE}_*^{i} \Big( \exp\big[\,t\, \sp{\L^{\tilde\xi}(t),r}\big]\,
= \widetilde{\EE}_*^{i} \Big( e^{t\, \sp{\L^{\tilde\xi}(t),r}}\,
I\{\s(\tilde\xi(t))=j\} \Big)\,,
\ee
which is a version of the \emph{Feynman-Kac formula}. Indeed, consider the function
$u(t,i)=\EE^{i} ( Z_j(t) )$ for fixed $j$. Since $u(t,i)=
(e^{t\A})_{ij}$,
it follows that $u(t,i)$ is the unique solution of the Cauchy problem
\[
\frac{d}{dt}\, u(t,i)=\sum_{k\in S} \tilde g_{ik}\, u(t,k) + r_i\, u(t,i), \quad
u(0,i) =\delta_{ij}\,,
\]
which is given by the Feynman-Kac formula.

(b) Summing over $j$ in (\ref{FC})
and using Varadhan's lemma of large deviation
theory (see \cite[p.~32]{deHo00} or \cite[Theorem 2.1]{Vara88}) together with
(\ref{popsize}) we arrive at the \emph{variational principle}
\[
\lambda=\lim_{t\ti}\frac1t\log \EE^i\big(|X(t)|\big)
= \max_{\nu\in \cP(S)}\big[ \sp{ \nu, r}- I_{\widetilde{\G}}(\nu)\big]\,,
\]
where $I_{\widetilde{\G}}$ is the large deviation rate function for the
empirical distribution of the Markov chain with transition rates $\tilde g_{ij}$;
cf.~(\ref{IG}) for its definition in the case of the transition rates $g_{ij}$.
In fact, it is not difficult to see that the maximum is attained at 
(and only at) the ancestral distribution $\a$.
This variational principle is behind the one found in \cite{HRWB02}.

(c) Just as in Remark \ref{conseq}(b) we find that
\[
\widetilde{W}(t) := \sum_{x \in X(t)} e^{-t \sp{\L^x(t),r}}
\]
is a martingale. In this martingale (which does
not seem to have been considered so far), each individual
at time $t$ is weighted according
to the mean fertility of its lineage.

\end{rem}

\section{Proofs}
\label{sec:proofs}
\setcounter{equation}{0}

\subsection{Transforming the tree}
\label{subsec:treetrans}

Here we prove Theorems \ref{size_bias} and \ref{KS}(b).
For the former we do not need that $\lambda$ is positive.

\begin{proof}[Proof of Theorem \ref{size_bias}]
It is sufficient to show that
 \be{size_bias_x}
 \hEE_*^{i}  \Big (F(\hX[0,t],\xi[0,t] )\,;\, \xi(t) = x \Big ) =
e^{-\lambda t}\,h_i\inv\,h_{\s(x)}  \;
\EE^{i} \Big ( F(X[0,t],x[0,t])\,;\, x \in X(t) \Big )
\ee
for all $x \in \XX$; the theorem then follows by summation over all
$x\in\XX$.
Suppose that $x = (i_1,\ldots,i_n; \l_1,\ldots,\l_n) \in \XX_n$, and
let $x_k=(i_1,\ldots,i_k; \l_1,\ldots,\l_k)$ be its ancestor
in generation $k$, $0\le k\le n$.

Consider the right-hand side of (\ref{size_bias_x}) and write
$e^{-\lambda t}\,h_i\inv\,h_{\s(x)}=q_1\,q_2\,q_3$ with
\[
q_1 =  e^{ - \lambda t}\,\prod_{k=0}^{n-1} c_{i_k}\,,\quad
q_2 = \prod_{k=0}^{n-1}\,\frac{\sp{ N_{x_k} ,h} }{c_{i_k}\,h_{i_k}}\,,\quad
q_3 = \prod_{k=0}^{n-1}\frac{h_{i_{k+1}}}{\sp{ N_{x_k} , h}} \,;
\]
of course, the random quantities $q_2$ and $q_3$ must then  be
included into the expectation. The factor $q_1$ corresponds to
the time change obtained when the exponential
parameter $a_{i_k}$ is replaced by $a_{i_k}{+}\lambda=a_{i_k}c_{i_k}$
along the ancestral line of $x$, i.e., when $\t_{x_k}$ is replaced by
$\htau_k$ for $k=0,\ldots,n$.
Indeed, the associated Radon-Nikodym density
is
\[
\bar{q}_1 = e^{ - \lambda T_x }\; \prod_{k=0}^{n} c_{i_k} \,.
\]
Conditioning $\bar{q}_1$ on the tree $X[0,t]$ up
to time $t$ and using the loss of memory property of the exponential law
of $\t_{x}$ we find that, almost surely on $\{T_{\tx} \leq t\}$,
\[
\EE^i\big(\,\bar{q}_1\,\big|\, X[0,t]\,\big)=
e^{ - \lambda t}\;  \EE^{i} \big( e^{-\lambda (T_{x}-t)} \, \big|\, T_{x}>t \big)
\; \prod_{k=0}^{n} c_{i_k} = q_1\,.
\]
Next, it is immediate from (\ref{sbd}) that the factor $q_2$ is precisely the
Radon-Nikodym density corresponding to a change from $N_{x_k}$
to the size-biased offspring $\hN_k$ for $k=0,\ldots,n-1$.
Finally, $q_3$ is the conditional selection probability for the trunk:
\[
q_3=\hPP^i\big(\xi_{k+1}=x_{k+1} \; \text{for} \;
    0 \leq k < n\,\big| \,\hX[0,t]\,\Big)\,.
\]
The right-hand side of (\ref{size_bias_x}) is therefore equal to
\[
\begin{split}
\hEE_*^{i} \Big (  F(\hX[0,t],x[0,t]) \,;\;  \hT_{\tx} \le t < \hT_{x},&\;
 \xi_{k+1}=x_{k+1}\; \text{for} \; 0 \le k < n \Big )\\
 &=   \hEE_*^{i} \Big ( F(\hX[0,t],\xi[0,t] )\,;\, \xi(t) = x \Big ) \,,
\end{split}
\]
as was to be shown.
\end{proof}

\bigskip
In the rest of this paper we assume that $\lambda>0$.

\begin{proof}[Proof of Theorem \ref{KS}(b)]
The basic observation is that the martingale $W(t)\, = \,$ $\sp{ Z(t),h}
\,e^{-\lambda t}$ considered in Remark \ref{conseq}(b) converges to a finite
limiting variable $W\ge 0$ $\PP^i$-almost surely for each $i$. When combined with
Theorem \ref{KS}(a) to be proved below, this implies the asserted
convergence result.
The essential part of the proof consists in showing that $W$ is nontrivial
if and only if condition (\ref{ZlogZ}) holds. There are two possible routes
to achieve this.

Either one can consider a discrete time skeleton $\d\NN$
and simply apply the discrete-time version of the Kesten-Stigum theorem.
For this one has to check that condition (\ref{ZlogZ}) holds
if and only if $\EE^i(Z_j(\d)\log Z_j(\d))<\infty$ for all $i,j\in S$,
which can be done.

Or, more naturally, one can use Theorem \ref{size_bias} to extend the
conceptual proof
 of Lyons et al.\ \cite{LPP95} and
Kurtz et al.\ \cite{KLPP97}
directly to continuous time. We spell out some details for the convenience
of the reader.
As in \cite{LPP95},
one observes first that $W$ is nontrivial if and only if $\hPP^i$ is
absolutely continuous with respect to $\PP^i$ (with Radon-Nikodym
density $W/h_i$), which is the case if ond only if
\be{hWleinfty}
\limsup_{t\ti}\hW(t) <\infty\quad\mbox{$\hPP^i$-almost surely;}
\ee
here we have put a hat on $W$ to stress the change of the underlying measure.

To check that (\ref{hWleinfty}) is equivalent to (\ref{ZlogZ})
one notices first that (\ref{ZlogZ}) is equivalent to
\be{hlogN}
\EE\big(\log \sp{\hN_i,h}\big)<\infty\quad\mbox{for all $i\in S$,}
\ee
by the properties of $\log$ and Eq.~(\ref{sbd}).
Next one observes that $\hX(t)\setminus\{\xi(t)\}$
is a branching process with immigration at the split times of the trunk
$\xi(t)$.
Specifically, let $\hT_{(n)}:=\hT_{\xi_n}$ be the $n$-th split time
and $\hN_{(n)}=\hN_{\xi_n}$ the $n$-th offspring  of the trunk. The $\hN_{(n)}$
are independent (conditionally on the trunk),
with distribution $\hp_{\s(\xi_n)}$.

Suppose first (\ref{hlogN}) fails, and pick any $j\in S$ with
$\EE\big(\log \sp{\hN_j,h}\big)=\infty$. Consider the subsequence
$(\hT_{(n_l)})_{l\ge 1}$ of split times of the trunk for which
$\s(\xi_{n_l})=j$.
Since the random variables $\log \sp{\hN_{(n_l)},h}$ are i.i.d.\ with
infinite mean, a standard Borel-Cantelli argument shows that
$\limsup_{l\ti} l\inv \log \sp{ \hN_{(n_l)},h}=\infty$ almost surely. 
On the other hand, 
$\limsup_{l\ti}\hT_{(n_l)}/l<\infty$ a.s.\ because the differences
$\hT_{(n_{l+1})}-\hT_{(n_l)}$ are i.i.d.\ with finite mean. This gives
$$\limsup_{l\ti}\hW(\hT_{(n_l)})\ge \limsup_{l\ti}\sp{\hN_{(n_l)},h}\,e^{-\lambda \hT_{(n_l)}}=\infty \quad\mbox{a.s.,}$$
so that (\ref{hWleinfty}) fails.

Conversely, suppose  (\ref{hlogN}) holds. As in Section
\ref{sec:bias_tree}, we consider the offspring $\hX_{n+1}^\sharp$
 of the trunk
created at time $\hT_{(n)}$ having type counting measure $\hN_{(n)}$.
We also introduce the $\s$-algebra $\cT$ generated by the trunk variables
$\{\hT_{(n)}, \hN_{(n)}:n\ge0\}$ and use a tilde to characterize
the trunk-reduced quantities obtained by removing the trunk individuals
from the population. Then for each $t>0$ we obtain, with the notation
(\ref{Zxt}),
\[
\begin{split}
\hEE^i_*\big(\widetilde{W}(t)\big|\cT\big)
&
=\sum_{n:\,\hT_{(n)}\le t} e^{-\lambda \hT_{(n)}}
\;\hEE^i_*\Big( \sum_{x\in \widetilde{X}_{n+1}^\sharp}
\sp{Z(x,t),h}\, e^{-\lambda (t-\hT_{(n)})}\Big|\,\cT\Big)
\\&
=\sum_{n:\,\hT_{(n)}\le t} e^{-\lambda \hT_{(n)}}
\sp{\widetilde{N}_{(n)},h}
\end{split}
\]
by the martingale property of $W(t)$ applied to the descendant trees
$X(x,\cdot)$. Now, (\ref{hlogN}) and a Borel-Cantelli argument imply that
%\[%\be{hlogZ}
%\hEE^i_*\big(\log \sp{ \hN_n,h}\big)<\infty
%\quad\mbox{for all $i\in S$ and $n\ge1$,}
%\]
$n\inv \log \sp{ \hN_{(n)},h}\to0$ almost surely. On the other hand,
$\liminf_{n\ti}\hT_{(n)}/n>0$ by the law of large numbers, whence
\[
\sum_{n\ge 0} e^{-\lambda \hT_{(n)}}\sp{\widetilde{N}_{(n)},h} <\infty \quad\mbox{a.s.}
\]
This means that, conditionally on $\cT$, $\widetilde{W}(t)$ is a 
submartingale with
bounded expectation,  which gives (\ref{hWleinfty}) by the submartingale
convergence theorem and finishes the proof of Theorem \ref{KS}(b).
The final identity $\{W>0\} = \Surv$ a.s. follows from the trivial
inclusion $\{W>0\} \subset \Surv$ and the well-known fact  that
$q_i=\PP^i(W=0)$ solves the equation $q_i=\EE\big(\prod_{j\in S}q_j^{N_{ij}}\big)$
which has the extinction probabilities as unique non-trivial solution
\cite[p.~205, Eq.~(25)]{AtNe72}.
\end{proof}

\subsection{Laws of large numbers for population averages}
\label{subsec:LLN}

In this section we are concerned with laws of large numbers for population
averages. We state a general such law for discrete time skeletons and
then use it to prove Theorems \ref{KS}(a) and \ref{pop_average}.
Recall from (\ref{Xxt}) that, for $t,u>0$ and $x\in X(t)$,
the path $X(x,[t,t{+}u])=\big( X(x,t{+}s)\big)_{0\le s\le u}$ describes
the subtree of $x$-descendants during the time interval $[t,t{+}u]$.

\begin{prop}\label{dLLN}
Let $\d,u>0$, $i,j\in S$, and $f: D\big([0,u],\fP(\XX)\big)\to\RR$
be a measurable function with existing mean $c_j=
\EE^j\big(f\circ X[0,u]\big)$. Then
\[
\lim_{n\ti} \frac1{Z_j(n\d)} \sum_{x\in X_j(n\d)}
f\circ X(x,[n\d,n\d+u]) = c_j
\quad\mbox{$\PP^i$-almost surely on $\Surv$.}
\]
\end{prop}

\begin{proof}
This result follows essentially from Lemmas 3 and 4 in \cite{KLPP97}.
Since this reference contains no proof of the former, we provide
a proof here for the sake of completeness.

We assume first that $\d$ is so large that $u<\d$ and $\rho:=\EE^j(Z_j(\d))>1$.
Such a $\d$ exists because $\lambda>0$ and $\A$ is irreducible.
Let $\cF_{n\d}$ denote the $\s$-algebra
generated by $X[0,n\d]$.  Since $u<\d$, for each $n\ge1$ the random variables
$\ph_{n,x}:= f\circ X(x,[n\d,n\d+u])$ with $x\in X_j(n\d)$ are
$\cF_{(n+1)\d}$-measurable and, conditionally on $\cF_{n\d}$, i.i.d.\
with mean $c_j$.  This implies that the sequence $(\ph_l)_{l\ge1}$
on $\Surv$ obtained by enumerating first $\{\ph_{1,x}:x\in X_j(\d)\}$ in some
order, then $\{\ph_{2,x}:x\in X_j(2\d)\}$ and so on, is still i.i.d.\
with mean $c_j$.  The strong law of large numbers therefore implies
that $\lim_{k\ti} (1/k) \sum_{l=1}^k\ph_l =c_j$ $\PP^i$-almost surely
on $\Surv$, and thus in particular that the subsequence
\[
A_n :=\frac{1}{\Psi_n}\sum_{l=1}^n\sum_{x\in X_j(l\d)}\ph_{l,x}
\]
converges to $c_j$ $\PP^i$-almost surely on $\Surv$ as $n\ti$;
here $\Psi_n=\sum_{l=1}^n \psi_l$ with $\psi_l=Z_j(l\d)$.

Next, the sequence $(\psi_l)_{l\ge1}$ dominates a single-type discrete-time
Galton-Watson process with mean $\rho>1$, and the latter
survives precisely on $\Surv$.  By Lemma 4 of \cite{KLPP97},
it follows that $\liminf_{l\ti}\psi_{l+1}/\psi_l \ge\rho$ almost surely
on $\Surv$.  This implies that
\[
\limsup_{n\ti}\Psi_{n-1}/\psi_n
= \limsup_{n\ti}\sum_{l=1}^{n-1}\psi_l/\psi_n
<\infty
\]
almost surely on $\Surv$.  As
\[
\frac1{\psi_n} \sum_{x\in X_j(n\d)}\ph_{n,x}
=A_n+ (A_n-A_{n-1})\,\Psi_{n-1}/\psi_n\,,
\]
the proposition follows in the case of large $\d$.

If $\d>0$ is arbitrary, we choose some $k\in\NN$ such that $\d':=k\d$
is so large as required above.  Let $0\le l<k$.  Applying the
preceding result to each of the subtrees $X(x,[l\d,\infty\ro)$ with
$x\in X(l\d)$ and averaging, we then find that
\[
\lim_{n\ti}\frac{1}{\psi_{nk+l}} \sum_{x\in X_j((nk+l)\d)}
\ph_{nk+l,x} =c_j
\]
$\PP^i$-almost surely on $\Surv$, and the proof is complete.
\end{proof}

\bigskip
A typical application of the preceding proposition is the following
corollary. Consider the $X_j(s)$-averaged type counting measure
\be{Cjus}
C_{j,u}(s)= \frac{1}{Z_j(s)} \sum_{x\in X_j(s)} Z(x,s{+}u)
\ee
at time $s{+}u$, where $Z(x,s{+}u)$ is defined by (\ref{Zxt}).
Proposition \ref{dLLN} then immediately implies
the following corollary.

\begin{cor}\label{discrete_forward_average}
For any $\d,u>0$ and $i,j\in S$,
\[
C_{j,u}(n\d)\toti \EE^j(Z(u))
\quad\mbox{$\PP^i$-almost surely on $\Surv$.}
\]
\end{cor}

To pass from a discrete time skeleton to continuous time we will use
the following continuity lemma which follows also from Proposition \ref{dLLN}.

%: Lem:continuity
\begin{lem}\label{continuity}
Given $\eps>0$, there exists some $\d>0$ such that for all
$i,j\in S$ and $k\in\NN$ one has
\be{Cont:1}
\limsup_{n\ti} \sup_{n\d\le s\le(n{+}1)\d}\;\frac{\|Z(s)\|}{\|Z(n\d)\|}
< 1+\eps\,,
\ee
\be{Cont:2}
\liminf_{n\ti}\; \inf_{n\d\le s\le(n{+}1)\d}\;\frac{Z_j(s)}{Z_j(n\d)}
> 1-\eps\,,
\ee
and
\be{Cont:3}
\liminf_{n\ti} \;\inf_{n\d\le s\le (n{+}1)\d}\;
\inf_{k\d\le u\le(k{+}1)\d}\;
\frac{\displaystyle\sum_{y\in X_j(s)}\|Z(y,s{+}u)\|}%\bigg/
{\displaystyle\sum_{y\in X_j(n\d)}\|Z(y,n\d{+}k\d)\|}
> 1-\eps
\ee
$\PP^i$-almost surely on $\Surv$.
\end{lem}

\begin{proof}
We begin by proving the upper bound (\ref{Cont:1}).
For $n\d\le s\le(n{+}1)\d$ we can write
\[
\|Z(s)\|= \sum_{x\in X(n\d)} |X(x,s)| \le
\sum_{x\in X(n\d)}  M(x,[n\d,(n{+}1)\d])\;,
\]
where $M(x,[n\d,(n{+}1)\d])=\max_{n\d\le s\le (n{+}1)\d}  |X(x,s)|$. Hence
\[
\sup_{n\d\le s\le(n{+}1)\d} \frac{\|Z(s)\|}{\|Z(n\d)\|} \le \max_{j\in S}
\frac{1}{Z_j(n\d)} \sum_{x\in X_j(n\d)}  M(x,[n\d,(n{+}1)\d])\;.
\]
By Proposition \ref{dLLN}, the last expression converges to
$m(\d):=\max_{j\in S}\EE^j\big(M(0,[0,\d])\big)$
almost surely on $\Surv$.  Now, $M(0,[0,\d])$ is dominated by the total
size at time $\d$ of the modified branching process for which the random
variables
$N_{x,\s(x)}$ in Section \ref{sec:model} are replaced by $N_{x,\s(x)}\vee 1$,
so that each individual has at least one offspring of its
own type. The latter process has a finite generator matrix, say $\A^+$. Hence
$m(\d)\le\max_j (e^{\d \A^+}\1)_j\to 1$ as $\d\to0$. This completes
the proof of
(\ref{Cont:1}).

Next we note that (\ref{Cont:2}) follows from (\ref{Cont:3})
by setting $u=k=0$. So it only remains to prove (\ref{Cont:3}).
Let $n\d\le s\le(n{+}1)\d$ and $k\d\le u\le(k{+}1)\d$.
Considering only those individuals $y\in X(s)$ already alive at time
$n\d$ and still alive at time $(n{+}1)\d$, and only those descendants
$z\in X(y,s{+}u)$ living during the whole period $[(n{+}k)\d,(n{+}k{+}2)\d]$,
we obtain the estimate
\[
\sum_{y\in X_j(s)}\|Z(y,s{+}u)\|\ge
\sum_{x\in X_j(n\d)}I\{\t_{x,n\d}>\d\}
\sum_{z\in X(x,(n{+}k)\d)} I\{\t_{z,(n{+}k)\d}>2\d\}\,.
\]
Here we write $\t_{x,t}=\inf\{u>0:\s(x)\not\in X(t{+}u)\}=
T_x-t$
for the remaining life time of $x\in X(t)$
after time $t$.
Proposition \ref{dLLN} therefore implies that the left-hand side
of (\ref{Cont:3}) is at least
\be{Cont3-ratio}
\EE^j\Big(I\{\t_{j,0}>\d\} \sum_{z\in X(k\d)} I\{\t_{z,k\d}>2\d\} \Big)
\Big/ \EE^j\big(| X(k\d)| \big)
\ee
$\PP^i$-almost surely on $\Surv$. By the Markov property, the numerator is equal to
\[
\EE^j\Big(I\{\t_{j,0}>\d\} \sum_{z\in X(k\d)} \exp[-2\d a_{\s(z)}] \Big)
\ge e^{-2\d a}\; \EE^j\big(I\{\t_{j,0}>\d\}\; |X(k\d)| \big)
\]
with $a=\max_i a_i$. The ratio in (\ref{Cont3-ratio}) is therefore
not smaller than $e^{-2\d a}\, (1-\eps_k)$, where
\[
\eps_k=\EE^j\Big(I\{\t_{j,0}\le\d\}\; |X(k\d)| \Big)
\Big/ \EE^j\big(| X(k\d)| \big)\,.
\]
For $k=0$ we have $\eps_0=1-e^{-\d a_j}$. For $k\ge 1$
we can use Theorem  \ref{size_bias} to obtain
\[
\eps_k=\hEE^j_*\big(I\{\t_{j,0}\le\d\}\; h_{\s(\xi(k\d))}\inv \big)
\big/\, \hEE^j_*\big( h_{\s(\xi(k\d))}\inv \big)
\le \frac{\max_i h_i}{\min_i h_i}\; (1-e^{-\d (a{+}\lambda)})\,.
\]
Hence, if $\d$ is sufficiently small then the ratio in (\ref{Cont3-ratio})
is larger than $1-\eps$.
\end{proof}

\bigskip
We are now ready for the proofs of Theorem \ref{KS}(a) and \ref{pop_average}.

\begin{proof}[Proof of Theorem \ref{KS}(a)]
Essentially we reproduce here the argument of \cite{KLPP97}.
Let $\eps>0$ be given and $\eps'>0$ be such that, for every
$\nu\in\cP(S)$,
$\|\nu-\pi\|<\eps$ whenever $\|a\nu-\pi\|<\eps'$ for some $a>0$.
Let $\d>0$ be so small as required in Lemma \ref{continuity}.
According to (\ref{PF}),
we can choose some $u\in\d\NN$ so large that
\[
\big\|\EE^j\big(Z(u)\,e^{-\lambda u}\big) - h_j\,\pi\big\| <\eps'\, \min_{i\in S}h_i
\]
for all $j\in S$.
Corollary \ref{discrete_forward_average} then implies that,
$\PP^i$-almost surely on $\Surv$,
\[
\big \|C_{j,u}(s)\, e^{-\lambda u}-h_j\,\pi\big\| <\eps' \, \min_{i\in S}h_i
\]
for all sufficiently large $s\in\d\NN$. Writing $\Pi(t)=Z(t)/\|Z(t)\|$ and
$a(t)=\frac{\| Z(t)\|\,e^{-\lambda u}}{\sp{ Z(t{-}u),h}}$ for $t>u$,
we conclude that
\[
\big\| a(t)\,\Pi(t)-\pi \big\|
\le \frac{1}{\sp{ Z(t{-}u),h}}\; \sum_{j\in S}Z_j(t{-}u)
\; \big \|C_{j,u}(t{-}u)\,e^{-\lambda u} -h_j\,\pi\big\|
< \eps'
\]
and therefore $\big\| \Pi(t)-\pi \big\| <\eps$
for all sufficiently large $t\in\d\NN$ a.s.\ on $\Surv$.
Finally, using (\ref{Cont:1}) and (\ref{Cont:2}) we find that
$\Pi_j(t) > (1-2\eps)\pi_j -\eps$ for all $j\in S$ and all sufficiently large
\emph{real} $t$, again a.s.\ on $\Surv$. Since $\eps$ was arbitrary and
$\Pi(t),\pi\in\cP(S)$, this gives the desired convergence result.
\end{proof}

%: pf 3.1
\begin{proof}[Proof of Theorem \ref{pop_average}]
Recall the definition (\ref{eq:ancestral-average}) of
$A^u(t) \in\cP(S)$,
the $X(t)$-average of the ancestral type distribution at time $t{-}u$,
and let $\a^u\in\cP(S)$ be given by its coordinates
$\a^u_j=\pi_j\, \EE^j(\|Z(u)\|)\,e^{-\lambda u}$.
Since $\a^u\to\a$ as $u\ti$ by (\ref{popsize}), it is sufficient
to show that
\be{eq:popav}
\PP^i\Big( \forall\, u>0: A^u(t)\toti[t] \a^u \,\Big|\,\Surv\Big) = 1.
\ee

Fix any $j\in S$, $u>0$ and $\d>0$. By Corollary \ref{discrete_forward_average},
\[
\big\|C_{j,u}(s)\big\| \to \EE^j\big(\|Z(u)\|\big)
\quad\mbox{as $s\ti$ through $\d\NN$}
\]
$\PP^i$-almost surely on $\Surv$. Combining this with Remark
\ref{rem:pop_average} and Theorem \ref{KS}(a) we obtain,
writing again $\Pi_j(s):=Z_j(s)/\|Z(s)\|$,
\[
\begin{split}
A^u_j(s+u) &=
 \frac{Z_j(s)\, \|C_{j,u}(s)\|}{|X(s+u)|} =
\frac{ \Pi_j(s)\, \|C_{j,u}(s)\|}{\sum_{k\in S}\Pi_k(s) \;\|C_{k,u}(s)\|}\\[1ex]
 &\toti[\d\NN\ni s]
\pi_j\, \EE^j\big(\|Z(u)\|\big)\Big/\EE^\pi\big(\|Z(u)\|\big) = \a^u_j
\end{split}
\]
$\PP^i$-almost surely on $\Surv$.

Next let $\eps>0$ be given and $\d>0$ be chosen according to
Lemma \ref{continuity}. Applying the above to $u=k\d$ with
arbitrary $k\in\NN$ and using (\ref{Cont:1}) and (\ref{Cont:3})
we find that
\[
\PP^i\Big( \forall\, u>0\ \forall \,j\in S :
\liminf_{t\ti}A^u_j(t) > (1-2\eps)\a^u_j   \,\Big|\,\Surv\Big) = 1\,,
\]
where the $u$-uniformity in (\ref{Cont:3}) allows us to bring the
$u$-quantifier inside of the probability.
%$A^u_j(t) > (1-2\eps)\a^u_j -\eps$ for all $j\in S$ and all
%sufficiently large \emph{real} $t$, still a.s.\ on $\Surv$.
This gives (\ref{eq:popav}) because $\eps$ is arbitrary and
$A^u(t)$ and $\a^u$ are probability measures on $S$.
\end{proof}

\subsection{Application of large deviation theory}
\label{subsec:LDP}

In this section we prove Theorems \ref{time_average} and \ref{mutationhistory}.
The main tools are the representation theorem \ref{size_bias} and the Donsker-Varadhan
large deviation principle for the empirical process of the retrospective mutation chain. In fact, these two ingredients together imply a
large deviation principle for the type histories as follows.
For every $\bnu\in\cP_\Theta(\S)$ let
\[
H_\G(\bnu)= \sup_{t>0} H\big(\bnu_{[0,t]};\bmu_{[0,t]}\big)\big/t
\]
be the process-level large deviation rate function for the retrospective mutation chain. In the above,
$\bnu_{[0,t]}$ and $\bmu_{[0,t]}$ are the restrictions of $\bnu$
and $\bmu$ to the time interval $[0,t]$, and
$H\big(\bnu_{[0,t]};\bmu_{[0,t]}\big)$ is their relative entropy.
See \cite[Eq.~(4.4.28)]{DeSt89}; alternative expressions
can be found in \cite[Theorem 4.4.38]{DeSt89} and
\cite[Theorems 7.3 and 7.4]{Vara88}.

\begin{thm}\label{LDP}
For the empirical type evolution process $R^x(t)$ as in $(\ref{Rt})$
we have, for $i\in S$ and closed $F\subset \cP_\Theta(\S)$
\[
\limsup_{t\ti} \frac1t \log \EE^i\Big(
\sum_{x\in X(t)} I\{R^x(t)\in F\}\Big) \le  \lambda -
\inf_{\bnu\in F} H_\G(\bnu)\;,
\]
while for open $G\subset \cP_\Theta(\S)$
\[
\liminf_{t\ti} \frac1t \log \EE^i\Big(
\sum_{x\in X(t)}I\{R^x(t)\in G\}\Big) \ge  \lambda -
\inf_{\bnu\in G} H_\G(\bnu)\;.
\]
Moreover, the function $H_\G$ is lower semicontinuous
with compact level sets and attains its minimum $0$ precisely at $\bmu$.
\end{thm}

\begin{proof}
In view of Theorem \ref{size_bias}, for every measurable
$C\subset \cP_\Theta(\S)$ we have
\[
 \EE^i\Big(\sum_{x\in X(t)} I\{R^x(t)\in C\}\Big)
= h_i\,e^{\lambda t}\,\hEE^i_*\Big(I\{R^\xi(t)\in C\}\,h_{\s(\xi(t))}\inv\Big)\;.
\]
Since $\max_i|\log h_i|<\infty$, the $h$'s can be ignored on the exponential
scale. The theorem thus follows from the Donsker-Varadhan large deviation
principle; see \cite[p.37, Theorem 7.8]{Vara88} or
\cite[Theorem 4.4.27]{DeSt89}, for example.
\end{proof}

\bigskip
There is a similar large deviation principle on the level of empirical distributions. For $\nu\in\cP(S)$ let
\be{IG}
  I_{\G}(\nu) = \sup_{ v \in \lo 0,\infty\ro^S}
  \Big[- \sum_{i\in S} \nu_i (\G v)_i/v_i  \Big] =
\inf_{\bnu\in\cP_\Theta(\S):\,\bnu_0 =\nu} H_\G(\bnu)
\ee
be the level-two rate function of the retrospective mutation chain;
here we write $\bnu_0$ for the time-zero marginal distribution of $\bnu$.
(For the second identity see \cite[p.37, Theorem 7.9]{Vara88}.)
Then the following statement holds.

\begin{cor}\label{LDP2}
For any $i\in S$ and closed $F\subset \cP(S)$,
\[
\limsup_{t\ti} \frac1t \log \EE^i\Big(
\sum_{x\in X(t)} I\{L^x(t)\in F\}\Big) \le  \lambda -
\inf_{\nu\in F} I_\G(\nu)\;,
\]
while for open $G\subset \cP(S)$
\[
\liminf_{t\ti} \frac1t \log \EE^i\Big(
\sum_{x\in X(t)}I\{L^x(t)\in G\}\Big) \ge  \lambda -
\inf_{\nu\in G} I_\G(\nu)\;.
\]
Moreover, the function $I_\G$ is continuous and strictly convex
and attains its minimum $0$ precisely at $\alpha$.
\end{cor}

\begin{proof} Simply replace the process-level large deviation principle for
the retrospective mutation chain by the one
for its empirical distributions.
The latter can either be deduced from the former by the contraction principle,
see \cite[Theorems 2.3\,\&\,7.9]{Vara88}, or be proved directly as in
\cite[Section IV.4]{deHo00}.
\end{proof}

\bigskip
We are now ready for the proofs of Theorems \ref{time_average} and \ref{mutationhistory}.

\begin{proof}[Proof of Theorem \ref{mutationhistory}]
Let $d$ be a metric for the weak topology on $\cP_\Theta(\S)$.
To be specific, we let $d_\S$ denote the Skorohod metric on $\S$ 
(defined in analogy to the one-sided case considered in \cite[p.~117, 
Eq. (5.2)]{EtKu86}),
and $d$ be the associated Prohorov metric on $\cP_\Theta(\S)$; see
\cite[p.~96, Eq. (1.1)]{EtKu86}. For any fixed $\eps>0$ we consider 
the set $C=\{\bnu\in\cP_\Theta(\S): d(\bnu,\bmu)\ge\eps\}$,
the complement of the open $\eps$-neigborhood of $\bmu$.
In view of Remark \ref{rem:mutationhistory} we need to show that
\[
\Gamma(t,C) := \frac{1}{|X(t)|}\sum_{x\in X(t)} I\{R^x(t)\in C\} \toti[t] 0
\]
$\PP^i$-almost surely on $\Surv$. In the first part of the proof
we will establish this convergence along a discrete time skeleton
$\d\NN$, where $\d>0$ is arbitrary.

Since $C$ is closed and $H_\G$ has compact level sets and attains its minimum
$0$ at $\bmu$ only, the infimum $c:=\inf_{\bnu\in C} H_\G(\bnu)$ is strictly
positive. We can therefore choose a constant $\lambda>\gamma>\lambda-c$.
We write
\[
\Gamma(t,C) = \Big(\frac{e^{\gamma t}}{|X(t)|}\Big)
\quad\Big( e^{-\gamma t}\sum_{x\in X(t)} I\{R^x(t)\in C\}\Big)
\]
and show that each factor tends to $0$ along $\d\NN$ a.s.\ on $\Surv$.
In view of Corollary \ref{discrete_forward_average} and Theorem \ref{KS}(a),
\[
\begin{split}
\frac{|X((n{+}1)\d)|}{|X(n\d)|} &=
\sum_{j\in S} \frac{Z_j(n\d)}{\|Z(n\d)\|} \quad
\frac{1}{Z_j(n\d)}\sum_{x\in X_j(n\d)} |X(x,(n{+}1)\d)|\\
&\toti \sum_{j\in S} \pi_j \,\EE^j(|X(\d)|) = e^{\lambda\d}
\quad\mbox{a.s. on $\Surv$.}
\end{split}
\]
Hence $n\inv\log |X(n\d)| \to \lambda\d$ and therefore
$e^{\gamma n\d}/|X(n\d)|\to0$ a.s. on $\Surv$.
On the other hand, using  Markov's inequality and
Theorem \ref{LDP} we obtain for any $a>0$
\[
\limsup_{n\ti} \frac{1}{n\d}\, \log
\PP^i\bigg( e^{-\gamma n\d}\sum_{x\in X(n\d)} I\{R^x(n\d)\in C\}\ge a\bigg)
\le \lambda-c-\gamma <0\,.
\]
The Borel-Cantelli lemma thus shows that also the second factor of $\Gamma(t,C)$
tends to $0$ a.s. as $t\ti$ through $\d\NN$. We therefore conclude
that $\lim_{n\ti} \Gamma(n\d,C)=0$ a.s. on $\Surv$.

To extend this result to the full convergence $t\ti$ along all reals
we pick some $0<\eps'<\eps$ and let $C'$ be defined in terms of $\eps'$
instead of $\eps$.
Also, let $A$ be an arbitrary closed set
in $\S$, $\eps^*=\eps-\eps'$, and $A^*=\big\{\s\in\S:d_\S(\s,A)<\eps^*\big\}$
the $\eps^*$-augmentation of $A$. Then for any
two time instants $s,t$ with $s\le t\le s+\d$ and
every $y\in X(t)$ we can write
\[
\begin{split}
R^y(t)(A) &\le\ \frac1t\int_0^s I_A(\th_u\s(y)_{\tper})\,du + \frac{\d}{t}\\[0.5ex]
&\le\ R^{y(s)}(s)(A^*)
%\frac1s\int_0^s I_{A^*}(\th_u\s(y(s))_{\tper[s]})\,du\\&\qquad
+ \frac1s\int_0^s I\big\{u:
d_\S(\th_u\s(y(s))_{\tper[s]},\th_u\s(y)_{\tper})\ge \eps^*
\big\}\,du
+ \frac{\d}{t}\,.
\end{split}
\]
By the locality of the Skorohod metric $d_\S$, there exists a constant $c=c(\eps^*)$
such that $d_\S(\th_u\s(y(s))_{\tper[s]},\th_u\s(y)_{\tper})< \eps^*$
whenever the interval $[-u,s{-}u]$ on which these functions agree contains $[-c,c]$.
The second term in the last sum is therefore at most $2c/s$, whence
\[
R^y(t)(A)\le R^{y(s)}(s)(A^*) +\eps^*
\]
for sufficiently large $s$. This means that
$d\big(R^y(t),R^{y(s)}(s)\big)<\eps^*$ and therefore
\[
\{ R^y(t)\in C\} \subset \{R^{y(s)}(s)\in C' \}
\]
when $s$ is large enough. For such $s$ we obtain
\[
\begin{split}
\Gamma(t,C)-\Gamma(s,C') &\le
\Big(\frac{1}{|X(t)|}-\frac{1}{|X(s)|}\Big)|X(t)|\\
&\qquad+ \frac{1}{|X(s)|} \sum_{x\in X(s)}
I \{R^x(s)\in C' \}\big (|X(x,t)|-1\big)\\
&\le 1- \inf_{s\le t\le s{+}\d} |X(t)|/|X(s)|\\
&\qquad + \frac{1}{|X(s)|} \sum_{x\in X(s)}
\big (M(x,[s,s+\d])-1\big)\;,
\end{split}
\]
where $M(x,[s,s{+}\d])=\max_{s\le t\le s{+}\d} |X(x,t)|$ as in the
proof of Lemma \ref{continuity}. Setting $s=n\d$, letting $n\ti$ and using
Theorem \ref{KS}(a) and Proposition \ref{dLLN} we see that the last
term converges to $\EE^\pi\big (M(0,[0,\d])-1\big)$ a.s. on
$\Surv$. According
to the proof of (\ref{Cont:1}), this limit can be made arbitrarily
small if $\d$ is chosen small enough. In combination with (\ref{Cont:2})
and the first part of this proof, this shows that
$\limsup_{t\ti} \Gamma(t,C)\le a$
for every $a>0$ almost surely on $\Surv$. The proof is thus
complete.
\end{proof}

\begin{proof}[Proof of Theorem \ref{time_average}]
There are two possible routes for the proof. One can either
repeat the argument above by simply replacing
Theorem \ref{LDP} by Corollary \ref{LDP2}. Or one notices that
$L^x(t)$ is the time-zero marginal of $R^x(t)$ and
that the marginal mapping $\bnu\to\bnu_0$ is continuous in the topologies
chosen.
The latter fact is used for the
derivation of the level-two large deviation principle
from that on the process level by means of the contraction principle; see
\cite[p.\ 34]{Vara88}.
\end{proof}

\bigskip

{\em Acknowledgement.} It is our pleasure to thank Nina Gantert,
Peter Jagers, G\"otz Kersting, and Anton Wakolbinger for helpful discussions,
and invaluable references to the branching literature.
Financial support from the German Research Council (DFG) and
the Erwin Schr\"odinger International Institute for Mathematical
Physics in Vienna is gratefully acknowledged.

%\bibliography{../eb}

\newcommand{\noopsort}[1]{} \newcommand{\printfirst}[2]{#1}
  \newcommand{\singleletter}[1]{#1} \newcommand{\switchargs}[2]{#2#1}

\end{document}